\documentclass{agtart_a}
\pdfoutput=1

\title{Homeomorphisms which are Dehn twists on the boundary}

\author{Darryl McCullough}
\givenname{Darryl}
\surname{McCullough}
\address{Department of Mathematics\\
University of Oklahoma\\\newline
Norman, Oklahoma 73019\\USA}
\email{dmccullough@math.ou.edu}
\urladdr{http://www.math.ou.edu/~dmccullough/}

\volumenumber{6}
\issuenumber{}
\publicationyear{2006}
\papernumber{47}
\startpage{1331}
\endpage{1340}

\doi{}
\MR{}
\Zbl{}

\keyword{3--manifold}
\keyword{boundary}
\keyword{Dehn twist}
\keyword{handlebody}
\keyword{compression body}
\subject{primary}{msc2000}{57M99}
\subject{secondary}{msc2000}{57R50}

\received{31 August 2005}
\revised{6 July 2006}
\accepted{20 July 2006}
\published{20 September 2006}
\publishedonline{20 September 2006}
\proposed{}
\seconded{}
\corresponding{}
\editor{}
\version{}

\arxivreference{math.GT/0604358}

\makeatletter
\def\cnewtheorem#1[#2]#3{\newtheorem{#1}{#3}[section]
\expandafter\let\csname c@#1\endcsname\c@theorem}
\makeatother

\newtheorem{theorem}{Theorem}
\cnewtheorem{corollary}[theorem]{Corollary}
\cnewtheorem{lemma}[theorem]{Lemma}

\begin{document}

\begin{asciiabstract}
A homeomorphism of a 3--manifold M is said to be Dehn twists on the
boundary when its restriction to the boundary of M is isotopic to the identity
on the complement of a collection of disjoint simple closed curves in
the boundary of M.  In this paper, we give various results about such
collections of curves and the associated homeomorphisms. In particular, if
M is compact, orientable, irreducible and the boundary of M is a single torus,
and M admits a homeomorphism which is a nontrivial Dehn twist on
the boundary of M, then M must be a solid torus.
\end{asciiabstract}

\begin{htmlabstract}
A homeomorphism of a 3&ndash;manifold M is said to be Dehn twists on the
boundary when its restriction to &part; M is isotopic to the identity
on the complement of a collection of disjoint simple closed curves in
&part; M.  In this paper, we give various results about such
collections of curves and the associated homeomorphisms. In particular, if
M is compact, orientable, irreducible and &part; M is a single torus,
and M admits a homeomorphism which is a nontrivial Dehn twist on
&part; M, then M must be a solid torus.
\end{htmlabstract}

\begin{abstract}
A homeomorphism of a $3$--manifold $M$ is said to be Dehn twists on the
boundary when its restriction to $\partial M$ is isotopic to the identity
on the complement of a collection of disjoint simple closed curves in
$\partial M$.  In this paper, we give various results about such
collections of curves and the associated homeomorphisms. In particular, if
$M$ is compact, orientable, irreducible and $\partial M$ is a single torus,
and $M$ admits a homeomorphism which is a nontrivial Dehn twist on
$\partial M$, then $M$ must be a solid torus.
\end{abstract}

\maketitle

\section*{Introduction}
\label{sec:intro}

A homeomorphism $h$ of a compact $3$--manifold $M$ is said to be Dehn twists
on the boundary when its restriction to $\partial M$ is isotopic to the
identity on the complement of a collection of disjoint simple closed curves
in $\partial M$. If this collection is nonempty, and the restricted
homeomorphism is not isotopic to the identity on the complement of any
proper subset of the collection, then we say that $h$ is Dehn twists
\textit{about\/} the collection.  The restriction of $h$ to $\partial M$ is
then isotopic to a composition of nontrivial Dehn twists about the curves,
where for us a Dehn twist may be a power of a ``single'' Dehn twist. Note
that the minimality condition implies that each curve in the collection is
essential in $\partial M$, and no two of them are isotopic in $\partial
M$. Our first main result gives strong restrictions on the collection of
curves.

\begin{theorem}
Let $M$ be a compact orientable $3$--manifold which admits a homeomorphism
which is Dehn twists on the boundary about the collection $C_1,\ldots\,$,
$C_n$ of simple closed curves in $\partial M$. Then for each~$i$, either
$C_i$ bounds a disk in $M$, or for some $j\neq i$, $C_i$ and $C_j$ cobound
an incompressible annulus in $M$.\par
\label{thm:mainthm}
\end{theorem}

Our second main result gives structural information about such
homeomorphisms. It refers to Dehn twists about disks and annuli in $M$,
whose definition is recalled in \fullref{sec:mainthm}.
\begin{theorem} 
Let $M$ be a compact orientable $3$--manifold which admits a homeomorphism
which is Dehn twists on the boundary about the collection $C_1,\ldots\,$,
$C_n$ of simple closed curves in $\partial M$. Then there exists a
collection of disjoint imbedded disks and annuli in $M$, each of whose
boundary circles is isotopic to one of the $C_i$, for which some
composition of Dehn twists about these disks and annuli is isotopic to $h$
on $\partial M$.
\label{thm:Dehn twists}
\end{theorem}
That is, $h$ must arise in the most obvious way, by composition
of Dehn twists about a collection of disjoint annuli and disks with a
homeomorphism that is the identity on the boundary.

\fullref{thm:mainthm} and \fullref{thm:Dehn twists} yield strong
statements for specific classes of manifolds. For the case when $M$ is a
compression body, examined in \fullref{sec:cbodies}, a
homeomorphism which is Dehn twists on the boundary is actually isotopic to
a product of Dehn twists about disjoint annuli and disks.  This appears
in Oertel \cite{Oertel} for the case when $M$ is a handlebody. Another application
is the following:
\begin{corollary} Let $M$ be a compact orientable irreducible $3$--manifold
with boundary a torus. If $M$ admits a homeomorphism which is a Dehn twist
on $\partial M$, then $M$ is a solid torus and the homeomorphism is isotopic
to a Dehn twist about a meridian disk.
\label{coro:torus boundary}
\end{corollary}
In particular, the only knot complement in $S^3$ (indeed, the
only irreducible complement of a knot in any closed orientable
$3$--manifold) admitting a homeomorphism which is a nontrivial Dehn twist on
the boundary is the trivial knot.
\begin{proof}[Proof of \fullref{coro:torus boundary}]
By \fullref{thm:mainthm}, $C$ bounds a disk in $M$. Since $M$ is
irreducible, this implies that $M$ is a solid torus.  By
\fullref{thm:Dehn twists}, there is a Dehn twist about the meridian
disk which is isotopic on $\partial M$ to the original homeomorphism, and
since any homeomorphism which is the identity on the boundary of a solid
torus is isotopic to the identity, this Dehn twist and the original
homeomorphism must be isotopic.
\end{proof}

It appears that most of our results could be extended to the nonorientable
case, adding the possibility of Dehn twists about M\"obius bands in
\fullref{thm:mainthm}, \fullref{thm:Dehn twists} and 
\fullref{coro:comp_body}, but the proof of \fullref{lem:reducible}
would require the more elaborate machinery of uniform homeomorphisms, found
in McCullough \cite{Montreal} or Chapter~12 of Canary and McCullough \cite{CM} (in particular, Lemma~12.1.2
of \cite{CM} is a version of Lemma 1.4 of \cite{Montreal} that applies to
nonorientable $3$--manifolds). \fullref{coro:torus boundary} fails in
the nonorientable case, however. Not only can a nonorientable manifold with
torus boundary admit Dehn twists about M\"obius bands, but an annulus can
meet the torus boundary in such a way that a Dehn twist about the annulus
will be isotopic on the boundary torus to an even power of a simple Dehn
twist about one of its boundary circles.

Some of the work presented here is applied in the article \emph{Knot adjacency,
genus and essential tori\/} by E\,Kalfagianni and X-S\,Lin \cite{K-L}. We
are grateful to the authors of that paper for originally bringing the
possibility of results like \fullref{thm:mainthm} and \fullref{thm:Dehn
twists} to our attention.

The research in this article was supported in part by NSF grant
DMS-0102463.

\section[Proof of Theorem 1]{Proof of \fullref{thm:mainthm}}
\label{sec:mainthm}

Recall that Dehn twists in $3$--manifolds can be defined as follows.
Consider first a properly imbedded and $2$--sided disk or annulus $F$ in a
$3$--manifold $M$. Imbed the product $F\times [0,1]$ in $M$ so that $(F\times [0,1])\cap
\partial M=\partial F\times[0,1]$ and $F\times\{0\}=F$. Let $r_\theta$
rotate $F$ through an angle $\theta$ (that is, if $F$ is a disk, rotate
about the origin, and if it is an annulus $S^1\times [0,1]$, rotate in the
$S^1$--factor). Fixing some integer $n$, define $t\co M\to M$ by $t(x)=x$
for $x\notin F\times[0,1]$ and $t(z,s)=(r_{2\pi ns}(z),s)$ if $(z,s)\in
F\times [0,1]$. The restriction of $t$ to $\partial M$ is a Dehn twist
about each circle of $\partial F$. Dehn twists are defined similarly when
$F$ is a $2$--sphere or a two-sided projective plane, M\"obius band, torus,
or Klein bottle (for the case of tori, there are infinitely many
nonisotopic choices of an $S^1$--factor to define $r_\theta$). Since a
properly imbedded closed surface in $M$ is disjoint from the boundary, a
Dehn twist about a closed surface is the identity on $\partial M$.

The proof of \fullref{thm:mainthm} will use the following result on
Dehn twists about annuli in orientable $3$--manifolds.

\begin{lemma}
Let $A_1$ and $A_2$ be properly imbedded annuli in an orientable
$3$--manifold $M$, with common boundary consisting of the loops $C'$ and
$C''$. Let $N'$ and $N''$ be disjoint closed regular neighborhoods in
$\partial M$ of $C'$ and $C''$ respectively, and let $t_i$ be Dehn twists
about the $A_i$ whose restrictions to $\partial M$ are supported on $N'\cup
N''$. If the restrictions of $t_1$ and $t_2$ to $N'$ are isotopic relative
to $\partial N'$, then their restrictions to $N''$ are isotopic relative to
$\partial N''$. Consequently, if $A$ is a properly imbedded annulus whose
boundary circles are isotopic in $\partial M$ (in particular, if they are
contained in a torus boundary component of $M$), then any Dehn twist about
$A$ is isotopic to the identity on~$\partial M$.
\label{lem:twisted}
\end{lemma}

\begin{proof}
The result is clear if the $A_i$ have orientations so that their induced
orientations on $C'\cup C''$ are equal, since then the imbeddings of
$S^1\times I \times [0,1]$ into $M$ used to define the Dehn twists can be
chosen to agree on $S^1\times \partial I \times [0,1]$. So we assume that
the oriented boundary of $A_1$ is $C'\cup C''$ and the oriented boundary of
$A_2$ is $C'\cup (-C'')$.

By assumption, $t_1$ and $t_2$ restrict to the same Dehn twist near
$C'$. Their effects near $C''$ differ in that after cutting along $C''$,
the twisting of $C''$ occurs in opposite directions, but since $M$ is
orientable, they also differ in that this twisting is extended to collar
neighborhoods on opposite sides of $C''$ (that is, the imbeddings of
$S^1\times \partial I\times [0,1]$ used to define the Dehn twists fall on
the same side of $C'$ but on opposite sides of $C''$). Each of these
differences changes a Dehn twist about $C''$ to its inverse, so their
combined effect is to give isotopic Dehn twists near $C''$.

The last remark of the lemma follows by taking $A_1=A$ and $A_2$ to be an
annulus with $\partial A_2=\partial A_1$, with $A_2$ parallel into
$\partial M$. All Dehn twists about $A_2$ are isotopic to the identity on
$\partial M$, so the same is true for all Dehn twists about $A_1$. 
\end{proof}

We will also need a fact about homeomorphisms of reducible
$3$--man\-i\-folds, even in many of the cases when $M$ itself is
irreducible.

\begin{lemma}
Let $W=P\# Q$ be a connected sum of compact orientable $3$--manifolds, with
$P$ irreducible. Let $S$ be the sum $2$--sphere.  Suppose that $\partial P$
is nonempty and that $g\co W\to W$ is a homeomorphism which preserves a
component of $\partial P$. Then there is a homeomorphism $j\co W\to W$,
which is the identity on $\partial W$, such that $jg(S)=S$.
\label{lem:reducible}
\end{lemma}

\begin{proof}
Let $P_1\#\cdots \# P_r\#R_1\#\cdots \#R_s$ be a prime factorization of
$W$, where each $P_i$ is irreducible and each $R_j$ is $S^2\times S^1$.
Let $\Sigma$ be the result of removing from a $3$--sphere the interiors of
$r+2s$ disjoint $3$--balls $B_1,\ldots\,$, $B_r$, $D_1$, $E_1$,
$D_2,\ldots\,$, $E_s$. For $1\leq i\leq r$, let $P_i'$ be the result of
removing the interior of a small open $3$--ball $B_i'$ from $P_i$, and
regard $W$ as obtained from $\Sigma$ and the union of the $P_i'$ by
identifying each $\partial B_i$ with $\partial B_i'$ and each $\partial
D_j$ with $\partial E_j$.

In \cite{Montreal} and in Section~12.1 of~\cite{CM}, certain \emph{slide
homeomorphisms\/} of $W$ are constructed. These can be informally described
as cutting $W$ apart along a $\partial B_i$ or $\partial D_j$, filling in
one of the removed $3$--balls to obtain a manifold $Y$, performing an
isotopy that slides that ball around a loop in the interior of $Y$,
removing the $3$--ball and gluing back together to obtain a homeomorphism
of the original $W$. Slide homeomorphisms are assumed to be the identity on
$\partial W$ (this is ensured by requiring that the isotopy that slides the
$3$--ball around the loop in $Y$ be the identity on $\partial Y$ at all
times). Lemma~1.4 of \cite{Montreal}, essentially due to M\,Scharlemann,
says that if $T$ is a collection of disjoint imbedded $2$--spheres in the
interior of $W$, then there is a composition $j$ of slide homeomorphisms
such that $j(T)\subset \Sigma$.

Since $P$ is irreducible, we may choose notation so that $P=P_1$ and
$S=\partial B_1$. Applying Lemma~1.4 of \cite{Montreal} with $T=g(S)$,
we obtain $j$ so that $jg(S)\subset\Sigma$. In particular, there is a
component $Z$ of $W-jg(S)$ whose closure contains $P_1'$. Since $g$ is
assumed to preserve a component of $\partial P_1$, the closure of $Z$ must
be $jg(P_1')$. Since $P_1$ is irreducible, $jg(S)$ must be isotopic to $S$
in $W$, so changing $j$ by isotopy we obtain $jg(S)=S$.
\end{proof}

We can now prove \fullref{thm:mainthm}.  Let $N_j$ be disjoint closed
regular neighborhoods of the $C_j$ in $\partial M$, and let $F$ be the
closure of $\partial M - \cup_jN_j$.  By hypothesis, we may assume that $h$
is the identity on $F$. Let $M'$ be another copy of $M$, and identify $F$
with its copy $F'$ to form a manifold $W$ with boundary a union of tori,
each containing one $C_j$. Denote by $T_j$ the one containing $C_j$. Let
$g\co W\to W$ be $h$ on $M$ and the identity map on $M'$, so that on
each $T_j$, $g$ restricts to a nontrivial Dehn twist about~$C_j$.

Fix any $C_i$, and for notational convenience call it $C_1$. Put $W_1=W$ if $W$ is
irreducible. Otherwise, write $W$ as $W_1\# W_2$ where $W_1$
is irreducible and $T_1\subseteq \partial W_1$, and let $S$ be the sum
sphere. By \fullref{lem:reducible}, there is a homeomorphism $j$ of $W$
that is the identity on $\partial W$, such that $jg(S)=S$. Split $W$ along
$S$, fill in one of the resulting $2$--sphere boundary components to obtain
$W_1$ and extend $jg$ to that ball. This produces a homeomorphism $g_1$ of
$W_1$ that restricts on each boundary torus of~$W_1$ to a nontrivial Dehn
twist about one of the~$C_j$.

Assume first that $W_1$ has compressible boundary. Since $W_1$ is
irreducible, it is a solid torus with boundary $T_1$. The only nontrivial
Dehn twists on $T_1$ that extend to $W_1$ are Dehn twists about a meridian
circle, showing that $C_1$ bounds a disk in $W_1$, and hence a disk $E$
in~$W$. Since $C_1$ does not meet $F$, we may assume that $E$ meets $F$
transversely in a collection of disjoint circles. The intersection $X'$ of
$E$ with $M'$ has a mirror image $X$ in $M$. Change $E$ by replacing $X'$
with $X$, producing a singular disk in $M$ with boundary $C_1$.  By the
Loop Theorem, $C_1$ bounds an imbedded disk in~$M$.

We call the argument in the previous paragraph that started with $E$ in $W$
and obtained a singular version of $E$ in $M$, having the same boundary as
the original $E$, a \emph{swapping\/} argument (since we are swapping pieces
of the surface on one side of $F$ for pieces on the other side).

Suppose now that $W_1$ has incompressible boundary. Let $V_1$ be
Johannson's characteristic submanifold of $W_1$ (\cite{Johannson}, also see
Chapter~2 of \cite{CM} for an exposition of Johannson's theory). Since
$\partial W_1$ consists of tori, $V_1$ admits a Seifert fibering and
contains all of $\partial W_1$ (in Johannson's definition, a component of
$V_1$ can be just a collar neighborhood of a torus boundary
component). Each $C_j$ in $W_1$ is noncontractible in $T_j$, and $T_j$ is
incompressible in $W_1$, so $C_j$ is noncontractible in $W_1$. This implies
that $C_j$ is noncontractible in $W$, hence also in~$M$.

It suffices to prove that $C_1$ and some other $C_i$ cobound an imbedded
annulus $A$ in $W_1$ and hence in $W$. For then, a swapping argument
produces a singular annulus in $M$ cobounded by $C_1$ and $C_i$. Since
$C_1$ and $C_i$ are noncontractible, a direct application of the
Generalized Loop Theorem \cite{Waldhausen1} (see
\cite[p.~55]{Hempel}) produces an imbedded annulus in $M$ cobounded by $C_1$
and~$C_i$.

By Corollary~27.6 of \cite{Johannson}, the mapping class group of $W_1$
contains a subgroup of finite index generated by Dehn twists about
essential annuli and tori. So by raising $g_1$ to a power, we may assume
that it is a composition of such Dehn twists. The Dehn twists about tori do
not affect $\partial W_1$, so we may discard them to assume that $g_1$ is a
composition $t_1\cdots t_m$, where each $t_k$ is a Dehn twist about an
essential annulus $A_k$. By Corollary~10.10 of \cite{Johannson}, each
$A_k$ is isotopic into $V_1$. By Proposition~5.6 of \cite{Johannson},
we may further change each $A_k$ by isotopy to be either horizontal or
vertical with respect to the Seifert fibering of $V_1$.

Suppose first that some $A_k$ is horizontal. Then $V_1$ is either
$S^1\times S^1\times I$ or the twisted $I$--bundle over the Klein bottle (a
horizontal annulus projects by an orbifold covering map to the base
orbifold, and the orbifold Euler characteristic shows that the base
orbifold is either an annulus, a M\"obius band, or the disk with two
order--$2$ cone points, the latter two possibilities yielding the two
Seifert fiberings of the twisted $I$--bundle over the Klein bottle). In the
latter case, $\partial V_1=T_1$, so $W_1=V_1$ and therefore $\partial
W_1=T_1$. By \fullref{lem:twisted}, each $t_k$ is isotopic to the
identity on $T_1$, hence so is $g_1$, a contradiction. So $V_1=S^1\times
S^1\times I$.

Since $A_k$ is horizontal, it must meet both components of $\partial V_1$,
and we have $V_1=W_1$ and $\partial W_1=T_1\cup T_i$ for some $i$. Let
$A_0=C_1\times I\subset S^1\times S^1\times I$. For an appropriate Dehn
twist $t$ about $A_0$, $t^{-1}g$ is isotopic to the identity on
$T_1$. Using Lemma~3.5 of \cite{Waldhausen}, $t^{-1}g$ is isotopic to a
level-preserving homeomorphism of $W_1$, and hence to the identity. We
conclude that $g_1$ is isotopic to $t$, and consequently $C_1$ and $C_i$
cobound an annulus in $W_1$.

It remains to consider the case when all $A_k$ are vertical. In this case,
each $t_k$ restricts on $\partial W_1$ to Dehn twists about loops isotopic
to fibers, so each $C_j$ in $\partial W_1$ is isotopic to a fiber of the
Seifert fibering on $V_1$.

Let $V_1'$ be the component of $V_1$ that contains $C_1$. Suppose first
that $V_1'\cap \partial W_1=T_1$. Then each $A_k$ that meets $T_1$ has both
boundary circles in $T_1$, so \fullref{lem:twisted} implies that $g_1$ is
isotopic to the identity on $T_1$, a contradiction. So $V_1'$ contains
another $T_i$. Since $C_1$ and $C_i$ are isotopic to fibers, there is an
annulus in $V_1'$ with boundary $C_1\cup C_i$.

\section[Proof of Theorem 2]{Proof of \fullref{thm:Dehn twists}}

\fullref{thm:mainthm} provides a properly imbedded surface $S$ which is
either an imbedded disk with boundary $C_n$ or an incompressible annulus
with boundary $C_n$ and some other $C_i$. For some Dehn twist $t_n$ about
$S$, $t_n$ and $h$ are isotopic near $C_n$. The composition $t_n^{-1}h$ is
isotopic on $\partial M$ to a composition of Dehn twists about
$C_1,\ldots\,$, $C_{n-1}$ (some of them possibly trivial). Induction on $n$
produces a composition $t$ as in the theorem, except for the assertion that
the disks and annuli may be selected to be disjoint.

Let $D_1,\ldots\,$, $D_r$ and $A_1,\ldots\,$, $A_s$ be the disks and annuli
needed for the Dehn twists in $t$. We first work on the annuli.

We will say that a union $\mathcal{A}$ of disjoint incompressible imbedded
annuli in $M$ is \emph{sufficient\/} for $A_1,\ldots\,$, $A_k$ if each
boundary circle of $\mathcal{A}$ is isotopic in $\partial M$ to a boundary
circle of one of the $A_i$, and if for any composition of Dehn twists about the set
$(\cup_{i=1}^rD_i)\cup(\cup_{j=1}^kA_j)$, there is a composition of Dehn
twists about the union of $\cup_{i=1}^rD_i$ and the annuli of $\mathcal{A}$
which has the same effect, up to isotopy, on $\partial M$. In particular,
${\mathcal{A}}=A_1$ is sufficient for $A_1$ alone. Inductively, suppose
that $\mathcal{A}$ is sufficient for $A_1,\ldots\,$, $A_{k-1}$.  By a
routine surgery process, we may change $\mathcal{A}$ so that $A_k$ and
$\mathcal{A}$ intersect only in circles essential in both $A_k$ and
$\mathcal{A}$.  (First, make $\mathcal{A}$ transverse to $A_k$. An
intersection circle which is contractible in $\mathcal{A}$ must also be
contractible in $A_k$, since both $\mathcal{A}$ and $A_k$ are
incompressible. If there is a contractible intersection circle, then there
is a disk $E$ in $A_k$ with $\partial E$ a component of $A_k\cap
\mathcal{A}$ and the interior of $E$ disjoint from $\mathcal{A}$. Replace
the disk in $\mathcal{A}$ bounded by $\partial E$ with $E$, and push off by
isotopy to achieve a reduction of $A_k\cap \mathcal{A}$.)

Now let $Z$ be a closed regular neighborhood of $A_k\cup
\mathcal{A}$. Since all intersection circles of $A_k$ with $\mathcal{A}$
are essential in both intersecting annuli, each component of $Z$ has a
structure as an $S^1$--bundle in which the boundary circles of
$\mathcal{A}$ and $A_k$ are fibers.

We will show that $Z$ contains a collection sufficient for $A_k\cup
\mathcal{A}$ and hence also for $A_1,\ldots\,$, $A_k$. We may assume that $Z$ is
connected. For notational simplicity, there is no harm in writing
$C_1,\ldots\,$, $C_m$ for the boundary circles of $\mathcal{A}$ and $A_k$,
since they are isotopic in $\partial M$ to some of the original $C_i$.

Fix a small annular neighborhood $N$ of $C_1$ in $Z\cap \partial M$.  Using
the $S^1$--bundle structure of $Z$, we can choose a collection
$B_2,\ldots\,$, $B_m$ of disjoint annuli, with $B_i$ running from $C_i$ to
a loop in $N$ parallel to $C_1$.

Consider one of the annuli $A$ of $A_k\cup \mathcal{A}$, say with boundary
circles isotopic to $C_i$ and $C_j$. If either $i$ or $j$ is $1$, say
$j=1$, then by \fullref{lem:twisted}, Dehn twists about $A$ have the same
effect on $\partial M$ as Dehn twists about $B_i$.  If neither is $1$, form
an annulus $B$ connecting $C_i$ to $C_j$ by taking the union of $B_i$,
$B_j$, and the annulus in $N$ connecting $B_i\cap N$ to $B_j\cap N$, then
pushing off of $N$ to obtain a properly imbedded annulus. Observe that any
Dehn twist about $B$ is isotopic on $M$ to a composition of Dehn twists
about $B_i$ and $B_j$. By \fullref{lem:twisted}, there is a Dehn twist
about $B$ whose effect on $\partial M$ is the same as the twist about
$A$. This shows that the collection $B_2,\ldots\,$, $B_m$ is sufficient for
$A_1,\ldots\,$, $A_k$ and completes the induction. So there is a collection
$\mathcal{A}$ sufficient for $A_1,\ldots\,$, $A_s$.

By further routine surgery, we may assume that each $D_i$ is disjoint from
$\mathcal{A}$. Then, surger $D_2$ to make $D_2$ disjoint from $D_1$, surger
$D_3$ to make it disjoint from $D_1\cup D_2$, and so on, eventually
achieving the desired collection of disjoint disks and annuli. 

\section{Compression bodies}
\label{sec:cbodies}

Compression bodies were developed by F\,Bonahon \cite{Bonahon}, in a
study of cobordism of surface homeomorphisms. They were used in work on
mapping class groups of $3$--manifolds \cite{MM,M,JDiff} and on deformations of hyperbolic structures on
$3$--manifolds \cite{CM}. The homeomorphisms of compression bodies
were further investigated by Oertel \cite{Oertel}, who develops an analogue
for compression bodies of the Nielsen--Thurston theory of surface
homeomorphisms. 

To fix notation and terminology, we recall that a \textit{compression body\/}
is a connected $3$--manifold $V$ constructed by starting with a compact
surface $G$ with no components that are $2$--spheres, forming $G\times
[0,1]$, and then attaching $1$--handles to $G\times\{1\}$. Compression
bodies are irreducible. They can be handlebodies (when no component of $G$
is closed) or product $I$--bundles (when there are no $1$--handles). The
\emph{exterior boundary\/} of $V$ is $\partial V-(G\times\{0\}\cup \partial
G\times [0,1))$. Note that if $F$ is the exterior boundary of $V$, and $N$
is a (small) regular neighborhood in $V$ of the union of $F$ with a
collection of cocore $2$--disks for the $1$--handles of $V$, then each
component of $\overline{V-N}$ is a product $X\times I$, where $X\times
\{0\}$ is a component of the frontier of $N$ and $X\times\{1\}$ is a
component of $G\times \{0\}$.

The following result was proven in \cite{Oertel} for the case of $V$ a
handlebody.

\begin{corollary}
Let $V$ be a compact orientable compression body, and let $h\co V\to V$
be a homeomorphism which is Dehn twists on the boundary about the
collection $C_1,\ldots\,$, $C_n$ of simple closed curves in $\partial V$.
Then $h$ is isotopic to a composition of Dehn twists about a collection of
disjoint disks and incompressible annuli in $V$, each of whose boundary
circles is isotopic in $\partial V$ to one of the~$C_i$.
\label{coro:comp_body}
\end{corollary}

To prove \fullref{coro:comp_body}, we note first that by
\fullref{thm:Dehn twists}, there is a composition $t$ of Dehn twists
about a collection of disjoint disks and incompressible annuli in $V$, such
that $t$ and $h$ are isotopic on $\partial V$. Changing $h$ by isotopy, we
may assume that $t^{-1}h$ is the identity on $\partial V$.
\fullref{coro:comp_body} is then immediate from the following lemma.

\begin{lemma}
Let $V$ be a compression body with exterior boundary $F$, and let $g\co
V\to V$ be a homeomorphism which is the identity on $F$. Then $g$ is
isotopic relative to $F$ to the identity.
\label{lem:identity}
\end{lemma}

\begin{proof}
We have noted that there is a collection of disjoint properly imbedded
disks $E_1,\ldots\,$, $E_n$, with boundaries in $F$, such that if $N$ is a
regular neighborhood of $F\cup(\cup_iE_i)$, then each component of
$\overline{V-N}$ is a product $X\times I$, where $X\times\{0\}$ is a
component of the frontier of $N$. Now $\partial E_1$ is fixed by $g$, so we
may assume that $g(E_1)\cap E_1$ consists of $\partial E_1$ and a
collection of transverse intersection circles. Since $V$ is irreducible, we
may change $g$ by isotopy relative to $F$ to eliminate these other
intersection circles, and finally to make $g$ fix $E_1$ as well as
$F$. Inductively, we may assume that $g$ is the identity on $F\cup
(\cup_iE_i)$ and then on $N$. Finally, for each component $X\times I$ of
$\overline{V-N}$, $g$ is the identity on $X\times\{0\}$. Using Lemma~3.5 of
\cite{Waldhausen}, $g$ may be assumed to preserve the levels $X\times\{s\}$
of $X\times I$, and then there is an obvious isotopy from $g$ to the
identity on $X\times I$, relative to $X\times\{0\}$. Applying these
isotopies on the complementary components of $N$, we make $g$ the identity
on $V$.
\end{proof}

\bibliographystyle{gtart}
\bibliography{link}

\end{document}